\documentclass{amsart}

\usepackage{amsmath}
\usepackage{amsfonts}
\usepackage{amssymb,enumerate}
\usepackage{amsthm}
\usepackage[all]{xy}

\newtheorem{lem}{Lemma}[section]
\newtheorem{cor}[lem]{Corollary}
\newtheorem{prop}[lem]{Proposition}
\newtheorem{thm}[lem]{Theorem}

\newtheorem{intthm}{Theorem}

\newtheorem{Defn}[lem]{Definition}
\newtheorem{Ex}[lem]{Example}
\newtheorem{Quest}[lem]{Question}
\newtheorem{Property}[lem]{Property}
\newtheorem{Properties}[lem]{Properties}
\newtheorem{Subprops}{}[lem]
\newtheorem{Para}[lem]{}
\newtheorem{Obs}[lem]{Observation}

\newtheorem{Remark}[lem]{Remark}

\newenvironment{defn}{\begin{Defn}\rm}{\end{Defn}}
\newenvironment{ex}{\begin{Ex}\rm}{\end{Ex}}
\newenvironment{quest}{\begin{Quest}\rm}{\end{Quest}}

\newenvironment{para}{\begin{Para}\rm}{\end{Para}}
\newenvironment{obs}{\begin{Obs}\rm}{\end{Obs}}
\newenvironment{remark}{\begin{Remark}\rm}{\end{Remark}}

\theoremstyle{definition}

\newcommand{\ideal}[1]{\mathfrak{#1}}
\newcommand{\m}{\ideal{m}}

\newcommand{\wcdim}{\operatorname{\PP_C\text{-}dim}}
\newcommand{\gcpd}{\operatorname{G_C\text{-}pd}}
\newcommand{\gc}{G_C}
\newcommand{\gpd}{\operatorname{G\text{-}pd}}

\newcommand{\pd}{\operatorname{pd}}

\newcommand{\id}{\operatorname{id}}

\newcommand{\im}{\operatorname{Im}}

\newcommand{\zz}{\mathbb{Z}}

\newcommand{\HH}{\operatorname{H}}
\newcommand{\Hom}{\operatorname{Hom}}

\newcommand{\coker}{\operatorname{Coker}}
\newcommand{\Ker}{\operatorname{Ker}}

\newcommand{\G}{\mathcal{G}}

\newcommand{\X}{\mathcal{X}}

\newcommand{\B}{\mathcal{B}}

\newcommand{\F}{\mathcal{F}}
\newcommand{\PP}{\mathcal{P}}

\newcommand{\tor}{\operatorname{Tor}}
\newcommand{\ext}{\operatorname{Ext}}

\newcommand{\gctor}{\operatorname{Tor}^{\G_C}}

\newcommand{\depth}{\operatorname{depth}}
\newcommand{\mcX}{\mathcal{X}}

\newcommand{\xra}{\xrightarrow}

\newcommand{\xpd}{\operatorname{\mcX\text{-}\pd}}
\newcommand{\pcpd}{\operatorname{\PP_C\text{-}\pd}}

\begin{document}
\dedicatory{Dedicated to the memory of Colleen Kilker}

\author{Diana White} \address{Department of Mathematics, LeConte College,
University of South Carolina, 1523 Greene St., Columbia, SC, 29208
USA}
   \email{dwhite@math.sc.edu}

\title[$G_C$-projective dimension]{Gorenstein projective dimension with respect to
a semidualizing module}


\begin{abstract}

We introduce and investigate the notion of $\gc$-projective modules
over (possibly non-noetherian) commutative rings, where $C$ is a
semidualizing module. This extends Holm and J{\o}rgensen's notion of
$C$-Gorenstein projective modules to the non-noetherian setting and
generalizes projective and Gorenstein projective modules within this
setting. We then study the resulting modules of finite
$\gc$-projective dimension, showing in particular that they admit
$\gc$-projective approximations, a generalization of the maximal
Cohen-Macaulay approximations of Auslander and Buchweitz. Over a
local ring, we provide necessary and sufficient conditions for a
$G_C$-approximation to be minimal.
% This paper contains crucial
%results for later papers that study and compare various relative
%(co)homology theories arising from semidualizing modules.
\end{abstract}
\maketitle

\section*{Introduction}\label{sec:intro}
Over a noetherian ring $R$, Foxby~\cite{foxby:gmarm},
Golod~\cite{golod:gdgpi}, and Vasconcelos~\cite{vasconcelos:dtmc}
independently initiated the study of semidualizing modules (under
different names): a module $C$ is semidualizing if $\Hom_R(C,C)\cong
R$ and $\ext_R^{\geqslant 1}(C,C)=0$.  Examples include the rank 1
free module and a dualizing (canonical) module, when one exists.
Golod~\cite{golod:gdgpi} used these to define $G_C$-dimension, a
refinement of projective dimension, for finitely generated modules.
The $G_C$-dimension of a finitely generated $R$-module $M$ is the
length of the shortest resolution of $M$ by so-called totally
$C$-reflexive modules; see Definition~\ref{totCref}. Motivated by
Enochs and Jenda's extensions in~\cite{enochs:gipm} of Auslander and
Bridger's G-dimension~\cite{auslander:smt}, Holm and
J{\o}rgensen~\cite{holm:sdmrghd} have extended this notion to
arbitrary modules over a noetherian ring.  The current paper
provides a unified and generalized treatment of these concepts, in
part by removing the noetherian hypothesis.  The tools developed in
this paper have been particularly useful for investigating the
similarities and differences between certain relative cohomology
theories~\cite{white:crct,white:gcac} and the stability properties
of operators on categories~\cite{white:sgc}.
%Some results in this
%paper, especially in Section 2, are (almost) immediate
%generalizations of the results by Holm in~\cite{holm:ghd}, which
%deals thoroughly with the case $C=R$. In these cases, we include the
%statements for completeness, but omit the proofs and refer the
%reader to~\cite{holm:ghd}.  We instead focus our efforts on dealing
%with the subtleties that often arise in the case $C\neq R$.  The
%tools developed in this paper are heavily used
%in~\cite{white:crct},~\cite{white:gcac}, and~\cite{white:sgc}.

Section~\ref{sec:gcproj} is devoted to the study of the
$G_C$-projective $R$-modules, which are built from projective and
$C$-projective modules; see Definition~\ref{cpltpcres}.  We show
that every module that is either projective or $C$-projective is
$G_C$-projective in Proposition~\ref{CR}.  In particular, every
$R$-module admits a $G_C$-projective resolution. Further properties
of the class of $G_C$-projective modules are contained in the
following result; see Theorem~\ref{projres}.

\begin{intthm}\label{introprojres}
The class of $G_C$-projectives is projectively resolving and closed
under direct summands. The class of finitely generated
$G_C$-projective $R$-modules is closed under summands. The set of
$G_C$-projective $R$-modules admitting a degreewise finite free
resolution is finite projectively resolving.
\end{intthm}

Section~\ref{sec:gcproj} ends with basic properties of the resulting
$G_C$-projective dimension.
%When proofs from~\cite{holm:ghd} follow
%through in a straight-forward manner, we omit them.
In particular, we show that, for an $R$-module $M$ of
$G_C$-projective dimension $n>0$, the $n$th kernel in any
$G_C$-projective resolution is $G_C$-projective.

Within the class of $G_C$-projective resolutions, the proper ones
exhibit particularly good lifting properties; see~\ref{proper}.
These are the subject of Section~\ref{sec:strict}.  Coupled with
Proposition~\ref{strictproper}, the following result shows that every
module of finite $G_C$-projective dimension admits a proper
$G_C$-projective resolution; see Theorem~\ref{strict}.

\begin{intthm}\label{introstrict} If $M$ is an $R$-module with finite
$G_C$-projective dimension, then $M$ admits a strict
$G_C$-projective resolution, that is, a $G_C$-resolution of the form
$$0 \to C\otimes_R P_n \to \cdots \to C\otimes_R P_1 \to G \to M \to
0$$ where $G$ is $G_C$-projective and $P_1,\ldots,P_n$ are
projective.
\end{intthm}

These strict $G_C$-projective resolutions give rise to
$G_C$-projective approximations, which are similar to the maximal
Cohen-Macaulay approximations of Auslander and Buchweitz
in~\cite{auslander:htmcma}.
%These proper resolutions give rise to
%an associated relative cohomology theory, see Remark~\ref{trivext2},
%which is studied in more depth
%in~\cite{white:crct},~\cite{white:gcac}.

Section~\ref{sec:fg} is concerned with comparing the
$G_C$-projective and totally $C$-reflexive properties; see
Definition~\ref{totCref}.  The next result is
Theorem~\ref{gctotCref}, which extends a result of Avramov,
Buchweitz, Martsinkovsky, and Reiten~\cite[(4.2.6)]{christensen:gd}.

\begin{intthm}\label{introgctotCref}
If $M$ and $\Hom_R(M,C)$ admit degreewise finite projective
resolutions, then $M$ is $G_C$-projective if and only if it is
totally $C$-reflexive.
\end{intthm}

The paper closes with several results on minimal proper
$G_C$-projective resolutions of finitely generated modules over
noetherian local rings.
%These minimal resolutions provide the
%crucial step in showing that various naturally arising cohomology
%theories are different, see~\cite{white:crct}.

\section{Preliminaries} \label{sec:back}

\emph{Throughout this work $R$ is a commutative ring with unity,
$\mcX=\mcX(R)$ is a class of unital $R$-modules, and $\mcX^f$ is the
subclass of finitely generated $R$-modules in $\X$.}

\

Homological dimensions built from resolutions are fundamental to
this work.  The prototypes are the projective and injective
dimensions.

\begin{para} \label{cxs}
An \emph{$R$-complex} is a sequence of $R$-module
homomorphisms
$$X = \cdots\xra{\partial^X_{n+1}}X_n\xra{\partial^X_n}
X_{n-1}\xra{\partial^X_{n-1}}\cdots$$
such that $ \partial^X_{n-1}\partial^X_{n}=0$ for each integer $n$; the
$n$th \emph{homology module} of $X$ is
$\HH_n(X)=\Ker(\partial^X_{n})/\im(\partial^X_{n+1})$.
A morphism of complexes $\alpha\colon X\to Y$ induces homomorphisms
$\HH_n(\alpha)\colon\HH_n(X)\to\HH_n(Y)$, and $\alpha$ is a
\emph{quasiisomorphism} when each $\HH_n(\alpha)$ is bijective.

The complex $X$ is \emph{bounded} if $X_n=0$ for $|n|\gg 0$; it is
\emph{acyclic} if $X_{-n}=0=\HH_n(X)$ for each $n>0$.   When $X$ is
acyclic, the natural morphism $X\to\HH_0(X)=M$ is a
quasiisomorphism, and $X$ is an \emph{$\mcX$-projective resolution}
of $M$ if each $X_n$ is in $\mcX$; in this event, the exact sequence
$$X^+ = \cdots\xra{\partial^X_{2}}X_1
\xra{\partial^X_{1}}X_0\to M\to 0$$ is the \emph{augmented
$\mcX$-projective resolution} of $M$ associated to $X$. Dually, one
defines $\X$-coresolutions and augmented $\X$-coresolutions. The
\emph{$\mcX$-projective dimension} of $M$ is defined as
$$\xpd_R(M)=\inf\{\sup\{n\mid X_n\neq 0\}\mid \text{$X$ is an
$\mcX$-projective resolution of $M$}\}.$$ The nonzero modules in
$\X$ are precisely the modules of $\xpd$ 0.
\end{para}

\begin{para} \label{projresolv}
The class $\mcX$ is \emph{projectively resolving} if
\begin{enumerate}[\quad(a)]
\item $\mcX$ contains every projective $R$-module, and
\item for every exact sequence of $R$-modules $0\to M'\to M\to M''\to 0$
with $M''\in\mcX$, one has $M\in\mcX$ if and only if $M'\in\mcX$.
\end{enumerate}

The class $\mcX$ is \emph{finite projectively resolving} if
\begin{enumerate}[\quad(a)]
\item $\mcX$ consists entirely of finitely generated $R$-modules,
\item $\mcX$ contains every finitely generated projective $R$-module, and
\item for every exact sequence of finitely generated $R$-modules $0\to M'\to M\to M''\to 0$
with $M''\in\mcX$, one has $M\in\mcX$ if and only if $M'\in\mcX$.
\end{enumerate}
\end{para}

\begin{para}\label{extensions}Consider an exact sequence of
$R$-modules $$0 \to M' \to M \to M'' \to 0.$$ The class $\X$ is
\emph{closed under extensions} when $M'$, $M'' \in\X$ implies
$M\in\X$, \emph{closed under kernels of epimorphisms} when $M$,
$M''\in \X$ implies $M'\in\X$ and \emph{closed under cokernels of
monomorphisms} when $M'$, $M\in\X$ implies $M'\in\X$.
\end{para}

\begin{para}\label{precover}
Let $M$ be an $R$-module. If $X\in\X$ and $\phi\colon X \to M$ is a
homomorphism, the pair $(X,\phi)$ is an $\X$-\emph{precover} of $M$
when, for every homomorphism $\psi\colon Y \to M$ where $Y\in\X$,
there exists a homomorphism $f\colon Y\to X$ such that $\phi
f=\psi$.  Enochs and Jenda introduced this terminology, which can be
found in~\cite{enochs:rha}.
\end{para}

\begin{para} \label{proper}
An $R$-complex $Z$ is \emph{$\mcX$-proper} if the complex
$\Hom_R(Y,Z)$ is exact for each $Y\in\mcX$.  If $\mcX$ contains $R$
and $Z$ is $\mcX$-proper, then $Z$ is exact.

An $\mcX$-resolution $X$ of $M$ is \emph{$\mcX$-proper} if the
augmented resolution
$X^+$ is $\mcX$-proper; by~\cite[(1.8)]{holm:ghd} $\mcX$-proper resolutions
are unique up to homotopy.  Accordingly, when $M$ admits an $\mcX$-proper
resolution $X$ and $N$ is an $R$-module, the
\textit{$n$th relative homology module} and the
\textit{$n$th relative cohomology module}
\begin{align*}
\tor_n^{\mcX}(M,N)&=\HH_n(X \otimes_R
N)&\ext^n_{\mcX}(M,N)&=\HH_{-n}\Hom_R(X,N)
\end{align*}
are well-defined for each integer $n$.
\end{para}

\begin{para}\label{fpinf}
A \emph{degreewise finite projective (respectively, free)
resolution} of an $R$-module $M$ is a projective (respectively,
free) resolution $P$ of $M$ such that each $P_i$ is a finitely
generated projective (respectively, free). Note that $M$ admits a
degreewise finite projective resolution if and only if it admits a
degreewise finite free resolution.  However, it is possible for a
module to admit a bounded degreewise finite projective resolution,
but not admit a bounded degreewise finite free resolution.  For
example, if $R=k_1\oplus k_2$, where $k_1$ and $k_2$ are fields,
then $M=k_1\oplus 0$ is a projective $R$ module, but it does not
admit a bounded free resolution.
\end{para}

The next result follows from well-known constructions, but the
author is unable to locate an elementary reference.

\begin{lem}\label{ffrclosed}
The class of $R$-modules admitting a degreewise finite projective
(respectively, free) resolution is closed under summands,
extensions, kernels of epimorphisms, and cokernels of
monomorphisms.\qed
\end{lem}

\begin{para} \label{sdm}
An $R$-module $C$ is \emph{semidualizing} if
\begin{enumerate}[\quad(a)]
\item $C$ admits a degreewise finite projective resolution,
\item The natural homothety map $R \to \Hom_R(C,C)$ is an isomorphism, and
\item $\ext^{\geqslant 1}_R(C,C)=0$.
\end{enumerate}
A free $R$-module of rank one is semidualizing. If $R$ is noetherian
and admits a dualizing module $D$, then $D$ is a semidualizing.

Note that this definition agrees with the established definition
when $R$ is noetherian, in which case condition (a) is equivalent to
$C$ being finitely generated. Also, since $\Hom_R(C,C)\cong R$ any
homomorphism $\phi\colon C^n \to  C^m$ can be represented by an $m
\times n$ matrix with entries in $R$.

Finally, note that the hypothesis that $C$ admits a degreewise
finite free resolution does not imply that $R$ is noetherian.  As
one example, take $R$ to be a non-noetherian ring and $C=R$.  For an
example with $C\neq R$, let $Q\to R$ be a flat local homomorphism of
commutative rings, with $Q$ noetherian and $R$ non-noetherian.  If
$C'$ is semidualizing over $Q$ with degreewise finite projective
resolution $F$, then $C=C'\otimes_Q R$ is semidualizing over the
non-noetherian ring $R$ with degreewise finite projective resolution
$F\otimes_Q R$.
\end{para}

\begin{para}\label{min}
Avramov and Martsinkovsky define a general notion of minimality for
complexes in~\cite[\S 1]{avramov:art}:  A complex $B$ is
\emph{minimal} if every homotopy equivalence $f\colon B\to B$ is an
isomorphism. Furthermore, by~\cite[(1.7)]{avramov:art} a complex $B$
is minimal if and only if every morphism $f\colon B\to B$ homotopic
to the identity map on $B$ is an isomorphism.
\end{para}

\begin{para}\label{tensorevaldefn}
Let $M$, $N$, and $F$ be $R$-modules.  The \emph{tensor evaluation}
homomorphism $$\omega_{MNF}\colon\Hom_R(M,N)\otimes_R F \to
\Hom_R(M,N\otimes_R F)$$ is defined by $\omega_{MNF}(\psi\otimes_R
f)(m)=\psi(m)\otimes_R f$.  It is straightforward to verify that
this is an isomorphism when $M$ is a finitely generated free (or
projective) $R$-module.
\end{para}

\begin{lem}\label{tensoreval}Let $F$ be a flat $R$-module.
\begin{enumerate}[\quad\rm(a)]
\item\label{t2} If $M$ admits a degreewise finite projective resolution $P$,
then for $i\geq 0$ there are isomorphisms $\ext^i_R(M,C\otimes_R
F)\cong\ext^i_R(M,C)\otimes_R F.$
\item\label{t3} If $M$ admits a degreewise finite projective resolution and
 $\ext^i_R(M,C)=0$ for some $i\geq 0$, then $\ext^i_R(M,C\otimes_R F)=0$.
\item\label{t4} If $M$ admits a degreewise finite projective resolution, $F$
is faithfully flat, and $\ext^i_R(M,C\otimes_R F)=0$ for some $i\geq
0$, then $\ext^i_R(M,C)=0$.
\end{enumerate}
\end{lem}

\begin{proof}

\eqref{t2} The maps $\omega_{P_iCF}$ are isomorphisms
by~\ref{tensorevaldefn}, hence the desired conclusion follows from
the flatness of $F$ and the resulting isomorphism of complexes
$$\Hom_R(P,C\otimes_R F)\cong\Hom_R(P,C)\otimes_R F.$$

\eqref{t3} and \eqref{t4}  These follow directly from \eqref{t2}.
\end{proof}

\begin{para}\label{cproj}
An $R$-module is \emph{$C$-projective} if it has the form
$C\otimes_R P$ for some projective $P$.  Set
$\PP_C=\PP_C(R)=\{C\otimes_R P \mid P \ \text{is projective}\}$.
These modules are studied extensively (in the non-commutative
setting) in~\cite{white:abc}.  We state for later use a Lemma that
follows readily from~\cite[(3.6, 5.6, 6.8)]{white:abc}.
\end{para}

%The following Lemma is contained in~\cite[(5.6)]{white:abc}.

\begin{lem}\label{wcr}Consider an exact sequence of
$R$-modules
\begin{align}\label{Cprojres}
0\to M' \to M \to M'' \to 0.
\end{align}
When $M''$ is a (finitely
generated) $C$-projective,  $M'$ is a (finitely generated)
$C$-projective if and only if $M$ is a (finitely generated)
$C$-projective. If all of the modules in~\eqref{Cprojres} are
$C$-projective, then~\eqref{Cprojres} splits.\qed
\end{lem}

\begin{para}\label{bass}The \textit{Bass class with respect to C},
 denoted $\B_C$ or $\B_C(R)$, consists of all $R$-modules $N$ satisfying
\begin{enumerate}[\quad\rm(a)]
\item\label{bass1} $\ext^{\geqslant 1}_R(C,N)=0$,
\item\label{bass2} $\tor_{\geqslant 1}^R(C,\Hom_R(C,N))=0$, and
\item\label{bass3} The evaluation map $\nu_{CN}\colon C\otimes_R\Hom_R(C,N)\to N$ is
an isomorphism.
\end{enumerate}
\end{para}

%\begin{proof}Suppose that $P$ and $P''$ are projective $R$-modules such
%that $M\cong C\otimes_R P$ and $M''\cong C\otimes_R P''$. By
%Lemma~\ref{wc}\eqref{wc1}, the modules $C\otimes_R P$ and
%$C\otimes_R P''$ are in $\B_C(R)$, and thus
%Lemma~\ref{bc}\eqref{bc5} implies that $M'$ is also in $\B_C(R)$,
%forcing the equality $\ext^1_R(C,M')=0$.

%Applying $\Hom_R(C,-)$ to the exact sequence~\eqref{Cprojres} and
%using~\ref{tensorevaldefn} together with the isomorphism
%$\Hom_R(C,C)\cong R$ yields an exact sequence
%\begin{align}\label{cpp}
%0\to\Hom_R(C,M') \to  P \to P'' \to 0
%\end{align}
%which splits as $P''$ is projective.  In particular,
%$P'=\Hom_R(C,M')$ is projective.
%Applying $C\otimes_R -$ to~\eqref{cpp} results in a split exact
%sequence
%\begin{align*}
%0\to C\otimes_R P' \to C\otimes_R P \to C\otimes_R P'' \to 0.
%\end{align*}
%As $M'$ is in $\B_C(R)$, one has $M'\cong
%C\otimes_R\Hom_R(C,M')\cong C\otimes_R P'$, which shows that $M'$ is
%$C$-projective and that the sequence~\eqref{Cprojres} splits.  In particular,
%this establishes the final statement of the theorem.

%Next, note that if $C\otimes_R P$ is a finitely generated
%$C$-projective, then since $M'\cong C\otimes_R P'$ is a summand of
%$C\otimes_R P$, it is also finitely generated.  The other
%implications are similar.
%\end{proof}

\section{$G_C$-projective modules}\label{sec:gcproj}

In this section we define and develop properties of $G_C$-projective
$R$-modules and the associated $G_C$-projective dimension. We begin
with a definition which extends the notion of $G_C$-projective
modules found in~\cite{holm:sdmrghd} (where they are referred to as
$C$-Gorenstein projective modules) to the non-noetherian setting.

\begin{defn}\label{cpltpcres} A \textit{complete $PC$-resolution}
is an exact sequence of $R$-modules
\begin{align}\label{cpltres}
X=\cdots \to P_1 \to P_0 \to C\otimes_R Q^0 \to C\otimes_R Q^1 \to
\cdots
\end{align}
where each $P_i$ and $Q^i$ is projective, and such that the complex
$\Hom_R(X,C\otimes_R Q)$ is exact for each projective $R$-module
$Q$.

An $R$-module $M$ is \textit{$G_C$-projective} if there exists a
complete $PC$-resolution as in~\eqref{cpltres} with $M\cong\coker(
P_1 \to P_0)$.
\end{defn}

Note that when $C=R$, the definitions above correspond to the
definitions of complete resolutions and Gorenstein projective
modules. The definition immediately gives rise to the following,
which generalizes~\cite[(2.3)]{holm:ghd}.

\begin{prop}\label{prop1}
A module $M$ is $G_C$-projective if and only if $\ext^{\geqslant
1}_R(M,C\otimes_R P)=0$ and $M$ admits a $\PP_C$-coresolution $Y$
with $\Hom_R(Y,C\otimes_R Q)$ exact for any projective $Q$.\qed
\end{prop}

\begin{obs}\label{obs1}If $M$ is a $G_C$-projective $R$-module, then
$M$ admits a complete $PC$-resolution of the form
\begin{equation}\label{cpltresfree}
\cdots \to F_1 \to F_0 \to C\otimes_R F^0 \to C\otimes_R F^1 \to
\cdots
\end{equation}
where each $F_i$ and $F^i$ is free.  To construct such a sequence
from a given complete $PC$-resolution argue as
in~\cite[(2.4)]{holm:ghd}.

When $X$ is a complex of the form~\eqref{cpltres}, then the complex
$\Hom_R(X,C\otimes_R Q)$ is exact for all projective $R$-modules $Q$
if and only if the complex $\Hom_R(X,C\otimes_R F)$ is exact for all
free $R$-modules $F$.  One implication is immediate.  For the other,
note that if $Q \oplus Q'$ is free, then we have the following
 isomorphism of complexes $\Hom_R(X,C\otimes_R (Q\oplus
Q'))\cong\Hom_R(X,C\otimes_R Q) \oplus \Hom_R(X,C\otimes_R Q')$.
\end{obs}

The next three results provide ways to create $G_C$-projective
modules.

\begin{prop}\label{cpltsums}
If $X_{\lambda}$ is a collection of complete $PC$-resolutions, then
$\coprod_{\lambda}X_{\lambda}$ is a complete $PC$-resolution. Thus,
the class of (finitely generated) $G_C$-projective $R$-modules is
closed under (finite) direct sums.
\end{prop}

\begin{proof}
For any projective $R$-module $Q$ there is an isomorphism,
$$\Hom_R\Bigl(\coprod_{\lambda}X_{\lambda},C\otimes_R Q\Bigr)\cong
\prod_{\lambda}\Hom_R(X_{\lambda},C\otimes_R Q).$$ Thus, if the
complex $\Hom_R(X_{\lambda},C\otimes_R Q)$ is exact for all
$\lambda$ then the complex
$\Hom_R\left(\coprod_{\lambda}X_{\lambda},C\otimes_R Q\right)$ is
exact.  It follows that a (finite) direct sum of (finitely generated)
$G_C$-projective $R$-modules is a (finitely generated) $G_C$-projective
$R$-module.
\end{proof}

\begin{lem}\label{exact}Let $P$ and $Q$ be projective $R$-modules and $X$ a
complex of $R$-modules.  If the complex $\Hom_R(X,C\otimes_R Q)$ is
exact, then the complex $\Hom_R(P\otimes_R X,C\otimes_R Q)$ is
exact. Thus, if $X$ is a complete $PC$-resolution of an $R$-module
$M$, then $P\otimes_R X$ is a complete $PC$-resolution of
$P\otimes_R M$. The converses hold when $P$ is faithfully projective.
\end{lem}

\begin{proof}
Assume the complex $\Hom_R(X,C\otimes_R Q)$ is exact.  Since
$\Hom_R(P,-)$ is an exact functor, the isomorphism of complexes
given by Hom-tensor adjointness
\begin{align*}
\Hom_R(P\otimes_R X,C\otimes_R Q) \cong\Hom_R(P,\Hom_R(X,C\otimes_R
Q))
\end{align*}
implies that $\Hom_R(P\otimes_R X,C\otimes_R Q)$ is exact.  It is
now straightforward to see that if $X$ is a complete $PC$-resolution
of an $R$-module $M$, then $P\otimes_R X$ is a complete
$PC$-resolution of $P\otimes_R M$.

If $P$ is faithfully projective, then the complex
$\Hom_R(P,\Hom_R(X,C\otimes_R Q))$ is exact if and only if the complex
$\Hom_R(X,C\otimes_R Q)$ is exact.
\end{proof}

\begin{prop}\label{CR}
If $P$ is $R$-projective, then $P$ and $C\otimes_R P$ are
$G_C$-projective.  Thus, every $R$-module admits a $G_C$-projective
resolution.
\end{prop}

\begin{proof}Using Lemma~\ref{exact}, it suffices to construct complete
$PC$-resolutions of $C$ and $R$.  By definition, $C$ admits an
augmented degreewise finite free resolution
$$X=\cdots \to R^{\beta_1} \to R^{\beta_0} \to C \to 0$$
and this is a complete $PC$-resolution of $C$. Indeed, the complex
$X$ is exact by definition and $C\cong\coker(R^{\beta_1}\to
R^{\beta_0})$. Furthermore, the complex $\Hom_R(X,C\otimes_R Q)$ is
exact for all projective $R$-modules $Q$ by
Lemma~\ref{tensoreval}\eqref{t3}, because $\ext_R^{\geqslant
1}(C,C)=0$. Thus, $C$ is $G_C$-projective.

We now show that $$\Hom_R(X,C)=0\to R \to C^{\beta_0} \to
C^{\beta_1} \to \cdots$$ is a complete $PC$-resolution of $R$.
First, left exactness of $\Hom_R(-,C)$ and the equality
$\ext^{\geqslant 1}_R(C,C)=0$ imply $\Hom_R(X,C)$ is exact.
Moreover, since $\Hom_R(X,C)$ consists of finitely presented
modules,
%Lemma~\ref{tensoreval}\eqref{t2} shows that,
for any projective $R$-module $Q$, tensor evaluation provides the
first isomorphism of complexes
$$\Hom_R(\Hom_R(X,C),C\otimes_R
Q) \cong\Hom_R(\Hom_R(X,C),C)\otimes_R Q \cong X\otimes_R Q.$$ The
second isomorphism follows from the fact that $\Hom_R(C,C)\cong R$.
These complexes are exact since the complex $X$ is exact and $Q$ is
flat.

Finally, since the class of $G_C$-projective $R$-modules contains
the class of projective $R$-modules, every $R$-module admits a
$G_C$-projective resolution.
\end{proof}

When $C=R$, the following proposition is contained
in~\cite[(2.3)]{holm:ghd}.  The proof is similar to that
of~\cite[(2.2)]{avramov:art}.

\begin{prop}\label{extvan}
If $X$ is a complete $PC$-resolution and $L$ is an $R$-module
admitting a bounded $\PP_C$-projective resolution, then the complex
$\Hom(X,L)$ is exact. Thus, if $M$ is $G_C$-projective, then
$\ext_R^{\geqslant 1}(M,L)=0$.\qed
\end{prop}

The following result is Theorem~\ref{introprojres} from the
introduction.

\begin{thm}\label{projres}
The class of $G_C$-projectives is projectively resolving and closed
under direct summands. The class of finite $G_C$-projective
$R$-modules is closed under summands.  The class of $G_C$-projective
$R$-modules admitting a degreewise finite projective resolution is
finite projectively resolving.
\end{thm}

\begin{proof}Consider an exact sequence
\begin{align}\label{ses}
0\to M' \xrightarrow{\iota} M \xrightarrow{\rho} M'' \to 0
\end{align}
 of $R$-modules. First,
assume that $M'$ and $M''$ are $G_C$-projective with complete
$PC$-resolutions $X'$ and $X''$, respectively. Use the Horseshoe
Lemmas in~\cite[(1.7)]{holm:ghd} and~\cite[(6.20)]{rotman:iha},
together with the fact that the classes of projective and
$C$-projective $R$-modules are closed under extensions to construct
a complex
$$X=\cdots\to P_1 \to P_0 \to C\otimes_R Q^0 \to C\otimes_R Q^1 \to
\cdots$$ with $P_i$ and $Q^i$ projective and a degreeswise split
exact sequence of complexes
$$0\to X' \to X \to X'' \to 0$$ such that $\coker(P_1\to P_0)\cong
M$. To show that $M$ is $G_C$-projective, it suffices to show that
$\Hom_R(X,C\otimes_R Q)$ is exact for all projective $R$-modules
$Q$. The sequence $$0 \to \Hom_R(X'',C\otimes_R Q) \to
\Hom_R(X,C\otimes_R Q) \to \Hom_R(X',C\otimes_R Q) \to 0$$ is an
exact sequence of complexes. Since the outer two complexes are
exact, the associated long exact sequence in homology shows that the
middle one is also exact.

Next, assume that $M$ and $M''$ are $G_C$-projective with complete
$PC$-resolutions $X$ and $X''$, respectively. Comparison lemmas for
resolutions, see e.g.~\cite[(1.8)]{holm:ghd} and
by~\cite[(6.9)]{rotman:iha}, provide a morphism of chain complexes
 $\phi\colon X \to X''$ inducing $\rho$ on the degree $0$ cokernels.
By adding complexes of the form $0 \to P''_i \xrightarrow{\id} P''_i
\to 0$ and $0 \to C\otimes_R (Q^i)'' \xrightarrow{\id} C\otimes_R
(Q^i)'' \to 0$ to $X$, one can assume $\phi$ is surjective. Since
both the class of projective and $C$-projective modules are closed
under kernels of epimorphisms, see Theorem~\ref{wcr}, the complex
$X'=\ker(\phi)$ has the form
$$X'=\cdots\to P'_1 \to P'_0 \to C\otimes_R (Q^0)' \to C\otimes_R
(Q^1)' \to \cdots$$ with $P'_i$ and $(Q^i)'$ projective. The exact
sequence $0\to X' \to X \to X'' \to 0$ is degreewise split by
Lemma~\ref{wcr}, so an argument similar to that of the previous
paragraph implies that $X'$ is a complete $PC$-resolution and $M'$
is $G_C$-projective.

Since the class of $G_C$-projective $R$-modules is projectively
resolving by the previous paragraphs and closed under arbitrary
direct sums by Proposition~\ref{cpltsums}, it follows from
Eilenberg's swindle~\cite[(1.4)]{holm:ghd} that they are also closed
under direct summands.

When the exact sequence~\eqref{ses} consists of modules admitting a
degreewise finite projective resolution, one can check that the
above constructions can be carried out using finite modules.
Finally, if $G$ is a finitely generated $G_C$-projective, then any
summand is also $G_C$-projective.  Since summands of finitely
generated modules are finitely generated, this implies that the
class of finitely generated $G_C$-projective modules is closed under
summands.
\end{proof}

When $C=R$, the next proposition follows readily from the symmetry
of the definition of the Gorenstein projectives.  However, the case
of $G_C$-projectives, the situation is more subtle. Nonetheless,
significant symmetry exists.

\begin{prop}\label{coker}
Every cokernel in a complete $PC$-resolution is $G_C$-projective.
\end{prop}

\begin{proof}Consider a complete $PC$-resolution
\begin{align}\label{7}
X=\cdots \to P_1 \to P_0 \to C\otimes_R Q^0 \to C\otimes_R Q^1 \to
\cdots
\end{align}
and set $M=\coker(P_1\to P_0)$ and $K=\coker(P_2\to P_1)$. Since $M$
and $P_0$ are $G_C$-projective, the exact sequence
\begin{align*}
0\to K \to P_0 \to M \to 0
\end{align*}
shows that $K$ is $G_C$-projective; see Theorem~\ref{projres}.
Inductively, one can show that $\coker(P_{i+1}\to P_i)$ is
$G_C$-projective for every positive integer $i$.

Set $N_{-1}=M$, $N_0=\coker(P_0\to C\otimes_R Q^0)$, and
$N_i=\coker(C\otimes_R Q^{i-1}\to C\otimes_R Q^i)$ for $i\geq 1$.
Using Proposition~\ref{prop1}, we will be done once we verify that
$\ext^{\geqslant 1}_R(N_i, C\otimes_R Q)=0$ for all projective
$R$-modules $Q$.  For each $i>-1$, consider the exact sequence
$$Y_i=0\to N_i\to C\otimes_R Q^{i+1} \to N_{i+1} \to 0.$$
By induction, one has $\ext^{\geqslant 1}_R(N_i, C\otimes_R Q)=0$.
Proposition~\ref{CR} implies that $C\otimes_R Q^{i+1}$ is
$G_C$-projective for each $i\geq 0$, and hence $\ext^{\geqslant
1}_R(C\otimes_R Q^{i+1}, C\otimes_R Q)=0$.  The long exact sequence
in $\ext_R(-,C\otimes_R Q)$ associated to $Y_i$ provides
$\ext^{\geqslant 2}_R(N_{i+1},C\otimes_R Q)=0$. Furthermore, since
$\Hom_R(X,C\otimes_R Q)$ is exact, so is the complex
$\Hom_R(Y_i,C\otimes_R Q)$.  Therefore, since $\ext^1_R(C\otimes_R
Q^{i+1},C\otimes_R Q)=0$, one has $\ext^1_R(N_{i+1}, C\otimes_R
Q)=0$.
\end{proof}

\label{sec:gcpd}

The class of $G_C$-projective $R$-modules can be used to define the
$G_C$-projective dimension, denoted $\gcpd_R(-)$; see~\ref{cxs}. The
following 5 results are proved similarly
to~\cite[(2.18),(2.19),(2.20),(2.21),(2.24)]{holm:ghd}.  We collect
them here for ease of reference.

\begin{prop}
Let $0 \to K \to G \to M \to 0$ be an exact sequence of $R$-modules
where $G$ is $G_C$-projective.  If $M$ is $G_C$-projective, then so
is $K$.  Otherwise, one has $\gcpd_R(K)=\gcpd_R(M)-1$.\qed
\end{prop}

\begin{prop}\label{gcpdses}
If $(M_{\lambda})_{\lambda\in\Lambda}$ is a collection of
$R$-modules, then
\begin{xxalignat}{3}&&{\hphantom{\square}} \gcpd_R\Bigl(\coprod_{\lambda}
M_{\lambda}\Bigr)=\sup\{\gcpd_R(M_{\lambda})\mid\lambda\in\Lambda\}.
&&\square
\end{xxalignat}
\end{prop}

\begin{prop}\label{gcpd}
Let $M$ be an $R$-module such that $\gcpd_R(M)$ is finite and let
$n$ be an integer.  The following are equivalent.
\begin{enumerate}[\quad\rm(i)]
\item $\gcpd_R(M)\leq n$.
\item $\ext^i_R(M,L)=0$ for all $i>n$ and all $R$-modules $L$ with
$\pcpd(L)<\infty$.
\item $\ext^i_R(M,C\otimes_R P)=0$ for all $i>n$ and all projective
 $R$-modules $P$.
\item In every exact sequence $0 \to K_n \to G_{n-1} \to \cdots \to
G_0 \to M \to 0$ where the $G_i$ are $G_C$-projective, one has that
$K_n$ is also $G_C$-projective.\qed
\end{enumerate}
\end{prop}
\begin{prop}\label{sup}
Let $M$ be an $R$-module with $\gcpd_R(M)<\infty$. If $M$ admits a
degreewise finite projective resolution, then there is an equality
$\gcpd_R(M)=\sup\{i\in \zz\mid\ext^i_R(M,C)\neq 0\}$.\qed

%If $M$ is an $R$-module of finite $G_C$-projective dimension and the
%modules $M$ and $\Hom_R(M,C)$ admit degreewise finite free
%resolutions, then $\gcpd_R(M)=\sup\{i\in \zz\mid\ext^i_R(M,C)\neq
%0\}$.\qed
\end{prop}

\begin{prop}\label{23}
If two modules in an exact sequence have finite $G_C$-projective
dimension, then so does the third.\qed
\end{prop}

When $C=R$, there are numerous proofs (see
e.g.~\cite[(3.4)]{avramov:art} or~\cite[(2.27)]{holm:ghd}) of the
following: if $M$ is an $R$-module of finite projective dimension,
then there is an equality $\pd_R(M)=\gpd_R(M)$.  Since
$G_C$-dimension can be viewed as a refinement of projective
dimension, it makes sense to ask the following:

\begin{quest}
If $M$ is an $R$-module of finite projective dimension, must
$\pd_R(M)=\gcpd_R(M)$?
\end{quest}

Over a noetherian, local ring, the affirmative answer in the case of
finitely generated modules follows immediately from the AB-formulas
for projective dimension and $G_C$-dimension.  Over a non-local
noetherian ring, an affirmative answer follows from work
in~\cite{holm:ghd} and~\cite{holm:sdmrghd}. However, as of the
writing of this paper, the author does not know the answer to this
question in general.

However, arguably the more natural comparison is between
$\PP_C$-dimension and $G_C$-dimension.  We have the following.

%We now compare the $\PP_C$-dimension, which is analogous
%to~\cite[(3.4)]{avramov:art}.

\begin{prop}
If $M$ is an $R$-module of finite $\PP_C$-projective dimension, then
$\pcpd_R(M)=\gcpd_R(M)$.
\end{prop}
\begin{proof}
Using Proposition~\ref{gcpd}, it suffices to show that if $M$ is
$G_C$-projective with finite $\PP_C$-projective dimension, then $M$
is $C$-projective.  To this end, consider an exact sequence of the
form $$0 \to K \to C\otimes_R P \to M \to 0$$ where $P$ is
projective and $\gcpd_R(K)<\infty$.  By Proposition~\ref{gcpd},
$\ext_R^1(M,K)=0$ so the above sequence splits, forcing $M$ to be a
summand of $C\otimes_R P$. Since the class of $C$-projectives is
closed under summands by~\ref{wcr}, this implies that $M$ is
$C$-projective, as desired.
\end{proof}

%\begin{remark}\label{trivext1} In~\cite{holm:sdmrghd}, Holm and
%J{\o}rgensen show that when $R$ is noetherian,
%$\gcpd_R(M)=\gdim_{R\ltimes C}(M)$, where $R\ltimes C$ is the
%trivial extension of $R$ by $C$. This provides a strong link between
%$G$-dimension and $G_C$-dimension, and in fact one can use this
%correspondence to prove most of the results from this section in the
%noetherian setting. At this time, we do not know if this
%correspondence holds when $R$ is not noetherian. See also
%Remark~\ref{trivext2}.
%\end{remark}

\section{$G_C$-projective resolutions and approximations}\label{sec:strict}

In this section we prove the existence of strict and proper
$G_C$-projective resolutions and of $G_C$-projective approximations.
These will give rise to well-defined relative (co)homology functors,
see Remark~\ref{trivext2}, which are further studied
in~\cite{white:crct} and~\cite{white:gcac}. We begin with the
requisite definitions.

\begin{defn}
Let $M$ be an $R$-module of finite $G_C$-projective dimension. A
\textit{strict $G_C$-projective resolution} of $M$ is a bounded
$G_C$-projective resolution $G$ such that for $i\geq 1$, there
exists a projective $R$-module $P_i$ such that $G_i\cong C\otimes_R
P_i$. This gives rise to an associated \emph{$G_C$-projective
approximation} of $M$; that is, an exact sequence of $R$-modules
$$0\to K\to G\to M \to 0$$ in which $\wcdim_R(K)$ is finite and $G$
is $G_C$-projective.
\end{defn}

We provide two examples.  The first corresponds to the situation
when $C$ is dualizing, the second to when $C=R$.

\begin{ex}
When $R$ is a local Cohen-Macaulay ring with dualizing module $D$,
Auslander and Buchweitz~\cite{auslander:htmcma} show that every
finitely generated module $M$ admits a maximal Cohen-Macaulay
approximation, that is, an exact sequence of the form
$$0\to K\to G\to M \to 0$$ where $K$ has finite injective dimension
and $G$ is maximal Cohen-Macaulay.  This gives rise to a resolution
of the form
$$0 \to D^{\alpha_n} \to \cdots\to
D^{\alpha_0} \to G \to M \to 0$$ where $G$ is a maximal
Cohen-Macaulay module.
\end{ex}

\begin{ex}
When $R$ is noetherian and $M$ is an $R$-module of finite
$G$-dimension, Avramov and Martsinkovsky~\cite[(3.8)]{avramov:art}
and Holm~\cite[(2.10)]{holm:ghd} provide several constructions of
$G$-approximations, that is, exact sequences of the form
$$0\to K\to G\to M \to 0$$ where $K$ has finite projective dimension and
$G$ is totally reflexive (see~\ref{totCref}). These give rise to
strict $G$-approximations, namely, exact sequences of the form
$$0 \to R^{\alpha_n} \to \cdots \to R^{\alpha_0} \to G \to M \to 0$$
where $G$ is totally reflexive.
\end{ex}

The existence of strict $G_C$-projective resolutions implies the
existence of proper $G_C$-projective resolutions.

\begin{prop}\label{strictproper}
Augmented strict $G_C$-projective resolutions are $\gc$-proper.
\end{prop}
\begin{proof}
Let $H$ be a $G_C$-projective $R$-module and
\begin{equation}\label{strictapprox}
0 \to C\otimes_R P_n \to \cdots\to  C\otimes_R P_1 \to G \to M \to 0
\end{equation}
an augmented strict $G_C$-projective resolution. Since $\ext^1_R(H,C
\otimes_R P_n)=0$ by Proposition~\ref{gcpd}, applying $\Hom_R(H,-)$
to the exact sequence $0\to C\otimes_R P_n \to C\otimes_R P_{n-1}
\to K_{n-2} \to 0$ provides an exact sequence
$$0\to \Hom_R(H,C\otimes_R P_n) \to \Hom_R(H,C\otimes_R P_{n-1}) \to
\Hom_R(H,K_{n-2})\to 0.$$ Continuing to break the exact
sequence~\eqref{strictapprox} into short exact sequences and
applying Proposition~\ref{gcpd} shows that~\eqref{strictapprox} is
$G_C$-proper.
\end{proof}

The existence of a strict $G_C$-projective resolution for a module
$M$ of finite $G_C$-projective dimension which is in the Bass class
of $R$ with respect to $C$ (see~\ref{bass}) was shown in
~\cite[(5.9)]{holm:sdmrghd}. We offer an alternative construction,
motivated by~\cite{auslander:htmcma}, that has the added advantage
of not requiring any Bass class assumption. When $R$ is noetherian
and $M$ is finitely generated, this is~\cite[(2.13)]{araya:hiasdb}.
We begin by proving a lemma.

\begin{lem}\label{pushout}
Let $\phi\colon G \to V$ be a homomorphism between $G_C$-projective
$R$-modules. If $0 \to G \xra{\psi} U \to N \to 0$ is an exact
sequence of $R$-modules such that $N$ is $G_C$-projective then the
pushout module $H$ of the maps $\phi$ and $\psi$ is
$G_C$-projective.
\end{lem}

\begin{proof}
We have a commutative diagram with exact rows
\[
\xymatrix{ 0 \ar[r]^{} & G \ar[d]^{\phi} \ar[r]^{\psi} & U \ar[d]^{}
\ar[r]^{}  &
N  \ar[d]^{=} \ar[r] & 0\\
0 \ar[r] & V  \ar[r]^{} & H \ar[r]^{} & N \ar[r] & 0.
\\}
\]
Since $N$ and $V$ are $G_C$-projective, Proposition~\ref{projres}
implies $H$ is $G_C$-projective.
\end{proof}

The next result contains Theorem~\ref{introstrict} from the
introduction.

\begin{thm}\label{strict} If $M$ is an $R$-module with finite $G_C$-projective
dimension, then $M$ admits a strict $G_C$-projective resolution and
hence a $G_C$-projective approximation.
\end{thm}

\begin{proof}
Assume $\gcpd_R(M)=n$.  By Proposition~\ref{gcpd}, truncating an
augmented free resolution of $M$ yields an augmented
$G_C$-projective resolution of $M$
$$0 \to G_n \xra{\phi_n} F_{n-1} \to \cdots \to F_0 \to M \to 0$$
A complete $PC$-resolution of $G_n$ gives rise to an exact sequence
$$0 \to G_n \xrightarrow{\psi} C\otimes_R P_n \to N \to 0$$ where
$P_n$ is projective and $N$ is $G_C$-projective. Lemma~\ref{pushout}
provides a commutative diagram (note that the orientation is not the
same as in the previous lemma)
\[
\xymatrix{ 0 \ar[r]^{} & G_n \ar[d]^{\psi} \ar[r]^{\phi_n} & F_{n-1}
\ar[d]^{} \ar[r]^{}  & F_{n-2} \ar[r] & \cdots \\
0 \ar[r] & C\otimes_R P_n  \ar[r]^{\phi'} & G_{n-1}  \\}
\]
with exact rows in which $G_{n-1}$ is $G_C$-projective.  As
$G_{n-1}$ is a pushout module, the maps $\phi_n$ and $\phi'$ have
isomorphic cokernels, resulting in a $G_C$-resolution
$$0 \to C\otimes_R P_n \xrightarrow{\phi '} G_{n-1} \to F_{n-2} \to \cdots \to F_0
\to 0.$$ Continuing this process yields a strict $G_C$-projective
resolution of $M$.
\end{proof}

\begin{remark}\label{trivext2}
As noted in the introduction, Proposition~\ref{strictproper} and
Theorem~\ref{strict} imply the following: every module $M$ of finite
$G_C$-projective dimension admits a proper $G_C$-projective
resolution. Hence, the relative (co)homology functors
$\ext_{G_C}^n(M,-)$ and $\gctor_n(M,-)$ are well-defined for each
integer $n$; see~\ref{proper}.
%In light of Remark 3.15, one may ask whether the entire theory of
%$\text{G}_C$-dimension reduces to G-dimension via the trivial
%extension $R\ltimes C$.  Even when $R$ is noetherian, this is not
%the case.  Indeed,
%This is where the limitations mentioned in Remark~\ref{trivext1}
%come into play, even over noetherian rings.  If $X$ is a proper
%$G_C$-resolution of a module $M$ over a noetherian ring $R$, then
%$X$ is a $G$-resolution of $M$ over the trivial extension $R\ltimes
%C$.  However, in general $X$ will not be a proper $G$-resolution,
%and so the study of the aforementioned relative (co)homology
%functors $\ext_{\G_C}^n(M,-)$ and $\gctor_n(M,-)$ does not reduce to
%the study of the the relative (co)homology functors
%$\ext_{\G}^n(M,-)$ and $\tor^{\G}_n(M,-)$.
\end{remark}

We close the section with a complement to Proposition~\ref{gcpdses},
which is proved as in~\cite[(2.11)]{holm:ghd}.

\begin{cor}
Let $0 \to G' \to G \to M \to 0$ be an exact sequence of
$R$-modules.  Assume $G$ and $G'$ are $G_C$-projective and that
$\ext^1_R(M,C\otimes_R Q)=0$ for all projective $R$-modules $Q$.
Then $M$ is $G_C$-projective.\qed
\end{cor}

\section{Connections with totally C-reflexive modules.}\label{sec:fg}

In this section, we reconnect with Golod's $G_C$-dimension.

\begin{defn}\label{totCref}Let $M$ be an $R$-module and assume that $M$
and $\Hom_R(M,C)$ admit a degreewise finite projective resolution.
The module $M$ is \emph{totally $C$-reflexive} if the following
conditions hold
\begin{enumerate}[\quad(a)]
\item The natural biduality map $M \to \Hom_R(\Hom_R(M,C),C)$ is an
isomorphism,
\item $\ext^{\geqslant 1}_R(M,C)=0$, and
\item $\ext^{\geqslant 1}_R(\Hom_R(M,C),C)=0$.
\end{enumerate}
\end{defn}

\begin{obs}
Finitely generated free modules are totally $C$-reflexive, as is the
$R$-module $C^n$ for any positive integer $n$. If $M$ is totally
$C$-reflexive, then it is straightforward to check that any summand
$M'$ of $M$ is also totally $C$-reflexive (using
Lemma~\ref{ffrclosed} to see that $M'$ admits a degreewise finite
free resolution).  Thus, finitely generated projective $R$-modules
are also totally $C$-reflexive, and so every finitely generated
$R$-module admits a resolution by totally $C$-reflexive modules.
\end{obs}

When $R$ is noetherian, the homological dimension which arises by
resolving a given module by totally $C$-reflexive modules is known
as the \emph{$G_C$-dimension} of a module, which was first
introduced by Golod; see~\cite{golod:gdgpi}.  In the case $C=R$,
this is Auslander and Bridger's
\emph{$G$-dimension}~\cite{auslander:smt}.

Next we provide a useful characterization of totally $C$-reflexive
modules, which generalizes~\cite[(4.1.4)]{christensen:gd}.

\begin{lem}\label{totCrefcplt}
An $R$-module $M$ is totally $C$-reflexive if and only if there is
an exact sequence of the form
\begin{equation}\label{cpltresfg2}
X=\cdots\to R^{\beta_1} \to R^{\beta_0} \to C^{\alpha_0} \to
C^{\alpha_1} \to \cdots
\end{equation}
with $M\cong\coker(R^{\beta_1}\to R^{\beta_0})$ and such that
$\Hom_R(X,C)$ is exact.
\end{lem}
\begin{proof}Set $(-)^{\dagger}=\Hom_R(-,C)$.
Assume first that $M$ is totally $C$-reflexive. By definition, there
exists augmented degreewise finite free resolutions
 $$F=\cdots\to R^{\beta_1} \to R^{\beta_0} \to M \to 0$$
$$G=\cdots\to R^{\alpha_1} \to R^{\alpha_0} \to M^{\dagger}\to 0.$$
% of $M$ and $\Hom_R(M,C)$, respectively.
The complexes $F^{\dagger}$ and $G^{\dagger}$ are exact, as
$\ext_R^{\geqslant 1}(M,C)=0=\ext_R^{\geqslant 1}(M^{\dagger},C)$.
The isomorphism $M\cong M^{\dagger\dagger}$  shows that
$G^{\dagger}$ has the form
$$G^{\dagger}\cong 0\to M \to C^{\alpha_0}\to C^{\alpha_1}\to \cdots.$$
Splicing together the complexes $F$ and $G^{\dagger}$ provides an
exact sequence $X$ of the form~\eqref{cpltresfg2} with
$M\cong\coker(R^{\beta_1}\to R^{\beta_0})$.  The fact that
$F^{\dagger}$ and $G$ are exact implies that $X^{\dagger}$ is exact.

Conversely, assume that $M$ admits a resolution $X$ of the
form~(\ref{cpltresfg2}) such that
\begin{align}\label{homxc}
X^{\dagger}=\cdots\to R^{\alpha_1} \to R^{\alpha_0} \to C^{\beta_0}
\to C^{\beta_1} \to \cdots
\end{align}
is exact.  Consider the following ``soft truncations'' of $X$
$$F=\cdots\to R^{\beta_1} \to R^{\beta_0} \to M \to 0$$
$$H=0\to M \to C^{\alpha_0} \to C^{\alpha_1} \to \cdots.$$
The complex $X^{\dagger}$ is exact and therefore so are
$F^{\dagger}$ and $H^{\dagger}$.

Since $F$ is an augmented free resolution of $M$, this implies that
$\ext_R^{\geqslant 1}(M,C)=0$.
%The exactness of $H^{\dagger}$ and the isomorphism
%$(C^{\alpha_i})^{\dagger}\cong R^{\alpha_i}$ imply that
%$H^{\dagger}$ is an augmented degreewise finite free resolution of
%$M^{\dagger}$.  Hence, the isomorphism $H\cong H^{\dagger\dagger}$
%implies that $\ext^{\geqslant 1}_R(M^{\dagger},C)=0$.
The biduality maps and exactness of $H^{\dagger}$ provide a
commutative diagram
\[
\xymatrix{ H=0 \ar[r]^{} & M \ar[d]^{\delta^C_M} \ar[r]^{} &
C^{\alpha_0} \ar[d]^{\delta^C_{C^{\alpha_0}}} \ar[r]^{}  &
C^{\alpha_1}
 \ar[d]^{\delta^C_{C^\alpha_1}} \ar[r] & \ldots \\
H^{\dagger\dagger} = 0 \ar[r] & \Hom_R(\Hom_R(M,C),C) \ar[r]^{} &
C^{\alpha_0} \ar[r]^{} & C^{\alpha_1} \ar[r] &\ldots.}
\]
The top row is exact by definition, while a routine diagram chase
and the fact that $\Hom_R(-,C)$ is left exact shows that the bottom
row is exact. Since $\delta^C_{C^{\alpha_1}}$ and
$\delta^C_{C^{\alpha_0}}$ are isomorphisms, the snake lemma implies
that the map $\delta^C_M$ is an isomorphism.  Finally, the exact
sequence $H^{\dagger}$ is an augmented degreewise finite free
resolution of $M^{\dagger}$.  Thus, exactness of
$H^{\dagger\dagger}$ implies that $\ext^{\geqslant
1}_R(M^{\dagger},C)=0$ and thus $M$ is totally $C$-reflexive.
\end{proof}

The next result is Theorem~\ref{introgctotCref} from the
introduction.

\begin{thm}\label{gctotCref}
If $M$ and $\Hom_R(M,C)$ admit degreewise finite projective
resolutions, then $M$ is $G_C$-projective if and only if it is
totally $C$-reflexive.
\end{thm}
\begin{proof}Set $(-)^{\dagger}=\Hom_R(-,C)$ and let $F$ and $G$ be
degreewise finite free resolutions of $M$ and $\Hom_R(M,C)$,
respectively.

Assume first that $M$ is totally $C$-reflexive.  By
Lemma~\ref{totCrefcplt}, there is an exact sequence
\begin{equation*}\label{cpltresfg}
X=\cdots\to R^{\beta_1} \to R^{\beta_0} \to C^{\alpha_0} \to
C^{\alpha_1} \to \cdots
\end{equation*}
with $M\cong\coker(R^{\beta_1}\to R^{\beta_0})$ and such that
$\Hom_R(X,C)$ is exact. An argument similar to the one used in the
proof of Lemma~\ref{tensoreval} implies that the complex
$\Hom_R(X,C\otimes_R P)$ is exact, and so $X$ is a complete
$PC$-resolution of $M$.

Conversely, assume that $M$ is $G_C$-projective and let $$\cdots \to
P_1 \to P_0 \to C\otimes_R F^0 \to C\otimes_R F^1 \to \cdots$$ be a
complete $PC$-resolution of $M$ in which each $F^i$ is a free
$R$-module; see Observation~\ref{obs1}. We show $M$ is totally
$C$-reflexive by constructing a complex $X$ as in
Lemma~\ref{totCrefcplt}.  To this end, it suffices to construct an
augmented $\PP_C^f$-coresolution
\begin{equation*}\label{cores}
Y=0\to M \to C\otimes_R R^{\alpha_0} \to C\otimes_R R^{\alpha_1} \to
\cdots
\end{equation*}
where each $\alpha_i$ is a non-negative integer and $Y^{\dagger}$ is
exact. Indeed, Proposition~\ref{gcpd} implies that $\ext^{\geqslant
1}_R(M,C\otimes_R P)=0$ for any projective $R$-module $P$; in
particular $\ext_R^{\geqslant 1}(M,C)=0$.  It follows that
$(F^+)^{\dagger}$ is exact.  Splicing together the complexes $F$ and
$Y$ provides the desired complex $X$.

We now build the complex $Y$ piece by piece.  Consider the exact
sequence $$0 \to M \to C\otimes_R F^0 \to G \to 0$$ arising from the
given complete $PC$-resolution of $M$. By Proposition~\ref{coker} we
know that $G$ is $G_C$-projective. Since $C\otimes_R F^0$ is a
direct sum of copies of $C$ we know that the image of the finitely
generated module $M$ is contained in a finite direct sum of copies
$C$.  That is, the image of $M$ is contained in a finitely generated
submodule $C\otimes_R R^{\alpha_0}$ of $C\otimes_R F^0$. Thus, we
have a commutative diagram with exact rows
\begin{equation}\label{sean}
\begin{split}
\xymatrix{ 0 \ar[r]^{} & M \ar[d]^{=} \ar[r]^{} & C\otimes_R
R^{\alpha_0}
\ar[d]^{} \ar[r]^{}  & H \ar[d]^{} \ar[r] & 0 \\
0 \ar[r] &  M \ar[r]^{} & C\otimes_R F^0 \ar[r]^{} & G \ar[r] &0. }
\end{split}
\end{equation}
Let $P$ be a projective $R$-module and set $\F=\Hom_R(-,C\otimes_R
P)$.  Since $C\otimes_R R^{\alpha_0}$ and $G$ are $G_C$-projective,
we have $\ext^1_R(G,C\otimes_R P)=0=\ext^1_R(C\otimes_R
R^{\alpha_0},C\otimes_R P)$.  Hence, applying $\F$ to~\eqref{sean}
yields a commutative diagram with exact rows
\[
\xymatrix{ 0 \ar[r]^{} & \F(G) \ar[d]^{} \ar[r]^{} & \F(C\otimes_R
F^0)
\ar[d]^{} \ar[r]^{}  & \F(M) \ar[d]^{=} \ar[r] & 0 \\
0 \ar[r] &  \F(H) \ar[r]^{} & \F(C\otimes_R R^{\alpha_0}) \ar[r]^{}
& \F(M) \ar[r] &\ext^1_R(H,C\otimes_R P) \ar[r]^{} & 0. }
\]
A routine diagram chase shows that $\ext^1_R(H,C\otimes_R P)=0$.
Proposition~\ref{gcpd} and Proposition~\ref{23} then imply that $H$
is $G_C$-projective. Since $M$ and $R^{\alpha_0}$ admit degreewise
finite projective resolutions, so does $H$. Applying $\Hom_R(-,C)$
to the exact sequence
$$0 \to M \to C\otimes_R R^{\alpha_0} \to H \to 0$$ gives rise to an
exact sequence $$0 \to \Hom_R(H,C) \to R^{\alpha_0} \to \Hom_R(M,C)
\to 0.$$ Here we used the facts that $\Hom_R(C\otimes_R
R^{\alpha_0},C)\cong R^{\alpha_0}$ and $\ext^1_R(H,C)=0$ because $H$
is $G_C$-projective. Since $\Hom_R(M,C)$ and $R^{\alpha_0}$ admit
degreewise finite projective resolutions, so does $\Hom_R(H,C)$; see
Lemma~\ref{ffrclosed}. Thus, we can proceed inductively to construct
the complex $Y$ with the given properties.
\end{proof}

\begin{cor}
If $M$ and $\Hom_R(M,C)$ admit degreewise finite projective
resolutions, then $M$ has finite $G_C$-projective dimension if and
only if it has finite $G_C$-dimension.  Moreover, these values
coincide.\qed
\end{cor}

Combining this with the AB-formula for $G_C$-dimension,
see~\cite[(3.14)]{christensen:sdctac}, and Proposition 2.16, we have
an AB-formula for modules of finite $\PP_C$-dimension.

\begin{cor}\label{ABformula}Let $R$ be a local, noetherian ring.
If $M$ is a finitely generated $R$-module of finite
$\PP_C$-dimension, then $\pcpd_R(M)=\depth(R)-\depth_R(M)$.\qed
\end{cor}

The next result compares with Theorem~\ref{strict}.

\begin{cor}\label{strictf}
If $M$ and $\Hom_R(M,C)$ admit degreewise finite projective
resolutions and $\gcpd_R(M)$ is finite, then $M$ admits a strict
$G_C^f$-resolution.\qed
\end{cor}

We conclude the paper with results on minimal proper
$G_C$-projective resolutions; see~\ref{min} for the definition of a
minimal complex. Note that Proposition~\ref{minimal}\eqref{diana}
shows, in particular, that such resolutions are strict.  We begin
with two lemmas, the first of which follows as
in~\cite[(8.1)]{avramov:art}.

\begin{lem}\label{cmin}
Over a local, ring $R$, a complex $H$ consisting of modules in
$\PP_C^f$ is minimal if and only if $\partial(H)\subseteq \m H$.\qed
\end{lem}

\begin{lem}\label{cminexists}
Let $R$ be local, noetherian and $M$ a finitely generated $R$-module
which admits a bounded $\PP_C^f$-resolution. Then $M$ admits a
minimal $\PP_C^f$-resolution.
\end{lem}

\begin{proof}
An augmented bounded $\PP_C^f$-resolution of $M$
\begin{align*}
X^+= \ 0 \to C^{\alpha_n} \to \ldots \to C^{\alpha_1} \to
C^{\alpha_0} \to M \to 0.
\end{align*}
is also an augmented strict $G_C$-projective resolution of $M$ and
so Proposition~\ref{strictproper} implies that it is proper.
Applying the functor $\Hom_R(C,-)$ to $X$ and using the fact that
$\Hom_R(C,C)\cong R$ yields an exact sequence
$$\Hom_R(C,X^+)= \  0 \to R^{\alpha_n} \to \ldots \to R^{\alpha_1} \to
R^{\alpha_0} \to \Hom_R(C,M) \to 0.$$
%Thus, $\Hom_R(C,X^+)$
which is an augmented finite free resolution of $\Hom_R(C,M)$. There
is an isomorphism of complexes $\Hom_R(C,X)\cong F\oplus G$ where
$F$ is an augmented minimal free resolution of $\Hom_R(C,M)$ and $G$
is a contractible complex of free modules.  Recall that $G$ is
\emph{contractible} if the identity map on $G$ is homotopic to the
zero map.

Since $M$ has finite $\PP_C$-dimension,~\cite[(2.9)]{white:cdim}
implies that $M\in\B_C(R)$ (see~\ref{bass} for the definition). This
provides the first isomorphism below
\begin{align*}
X^+ & \cong C\otimes_R \Hom_R(C,X^+)\\ &\cong (C\otimes_R F) \oplus
(C\otimes_R G)
\end{align*}
while the second follows from the isomorphism $\Hom_R(C,X^+)\cong
F\oplus G$ and the fact that finite direct sums commute with tensor
products. It is now straightforward to verify that the complex
$C\otimes_R F$ is contractible and that the complex $C\otimes_R F$
is a minimal $\PP_C$-resolution of $M$, as desired.
\end{proof}

The following structure result is the key to demonstating the
differences between the relative cohomology theories $\ext_{\PP_C}$,
$\ext_{G_C}$, and $\ext_R$ in~\cite{white:crct}.

\begin{prop} \label{minimal}
Assume that $R$ is local, and let $M$ be a finitely generated
$R$-module of finite $G_C$-projective dimension.  If $M$ and
$\Hom_R(M,C)$ admit degreewise finite projective resolutions, then
the following hold.
\begin{enumerate}[\quad\rm(a)]
\item The module $M$ admits a minimal proper $G_C$-projective
resolution.
\item\label{diana} A given $G_C$-projective resolution $H$ of $M$ is minimal if
and only if the following conditions hold.
\begin{enumerate}[\quad\rm(1)]
\item $H_n \cong C^{\alpha_n}$ for all $n\geq 1$,
\item $\partial^H_n(H_n)\subseteq\m H_{n-1}$ for all $n\geq 2$, and
\item\label{part3} $\partial^H_1(H_1)$ contains no nonzero $C$-summand of $H_0$.
\end{enumerate}
\end{enumerate}
\end{prop}

\begin{proof}
We begin by showing that a $G_C$-projective resolution $H$
satisfying conditions (1)--(3) is minimal.  First, observe that 
$H_0$ is finitely generated because $M$ and $H_1$ are so.  Let
$\gamma\colon H\to H$ be a morphism that is homotopic to $\id_H$. 
Using~\ref{min}, we need to show that $\gamma_n$ is an isomorphism
for each integer $n$.  

For $n\geq 0$, let  $\theta_n\colon H_n\to H_{n+1}$ be maps such
that
$\gamma_n-\id_{H_n}=\theta_{n-1}\partial_n^{H}+\partial_{n+1}^{H}\theta_n$,
which exist since $\gamma$ is homotopic to $\id_H$. For $n\geq 2$,
condition (2) implies $\partial_n^H\otimes_R k=0$, and so
$\gamma_n\otimes_R k-\id_{H_n}\otimes_R k=0$. Nakayama's Lemma
implies $\gamma_n$ is a surjective endomorphism, and hence
bijective.

Now let $n=1$.  We verify the containment
$\im(\theta_0\partial^H_1)\subseteq \m C^{\alpha_1}$ and then an
argument similar to that in the previous paragraph shows that
$\gamma_1$ is an isomorphism. Suppose
$\im(\theta_0\partial^H_1)\nsubseteq \m C^{\alpha_1}$. This means
the matrix representation of $\theta_0\partial_1^H$ contains a unit;
see~\ref{sdm}.  Thus, there exist maps $\rho\colon C^{\alpha_1}\to
C$ and $\iota\colon C \to C^{\alpha_1}$ such that
$\rho\theta_0\partial^H_1\iota=\id_C$. This provides a splitting
$\partial^H_1\iota$ of $\rho\theta_0\colon H_0\to C$, and so
$H_0\cong C \oplus\ker(\rho\theta_0)$.  Finally, the summand
$C\oplus 0$ is isomorphic to $\im(\partial^H_1\iota)$ which is
contained in $\im(\partial^H_1)\subseteq H_0$, contradicting
assumption (3).

The fact that $\gamma_0$ is an isomorphism now follows as
in~\cite[(8.5)]{avramov:art}.

Next we show that $M$ admits a resolution satisfying (1)--(3).  With
the first part of this proof, this will establish part (a). First,
note that by Corollary~\ref{strictf}, there exists a
$G_C$-projective approximation of $M$
\begin{align*}
Y=0\to K\to H \to M \to 0
%0\to C^{\alpha_n} \xrightarrow{\partial_n} \cdots
%\xrightarrow{\partial_1} C^{\alpha_1} \xrightarrow{\partial_0} H_0
%\to M \to 0
\end{align*}
where $H$ is totally $C$-reflexive
% a finitely generated $G_C$-projective
 and $K$ admits a bounded $\PP_C^f$-resolution.

%This gives rise to an exact sequence
%\begin{align}\label{strict3}
%Y'=0\to \im(\partial_1) \xrightarrow{\partial_0} H_0 \to M \to 0
%\end{align}
If possible, write $H\cong C\oplus H'$ and assume $K$ contains the
non-zero $C$-summand $C\oplus 0$ of $H$, say $K=C\oplus K'$ for some
$K'$.  One checks readily that the compatibility of the two
splittings gives rise to a split exact sequence of complexes,
written vertically
$$
\xymatrix{
0 \ar[d] & & 0 \ar[d] & 0 \ar[d] & 0 \ar[d] & \\
Y' \ar[d] & 0 \ar[r] & K' \ar[r] \ar[d] & H' \ar[r] \ar[d] & M \ar[r] \ar[d] & 0 \\
Y \ar[d] & 0 \ar[r] & K \ar[r] \ar[d] & H \ar[r] \ar[d]
& M \ar[r] \ar[d] & 0 \\
X \ar[d] & 0 \ar[r] & C \ar[r] \ar[d] & C \ar[r] \ar[d]
& 0  \\
0 & & 0 & 0. }$$ Since $Y'$ is a $G_C^f$-projective approximation of
$M$, one can repeat this process. Finitely many iterations yield a
$G_C$-approximation $Y_0=0 \to K_0 \to H_0 \to M \to 0$ where $K_0$
does not contain a nonzero $C$-summand of $H_0$.
Lemma~\ref{cminexists} implies that $K$ admits a minimal
$\PP_C^f$-resolution $Z$. Splicing together $Y_0$ and $Z$ at $K$
provides a resolution of $M$ satisfying (1)--(3).

Finally, let $G$ be a resolution of $M$ satisfying conditions
(1)--(3), and let $H$ be a minimal proper $G_C$-projective
resolution of M.  It follows from [11, (1.8)] that $G$ and $H$ are
homotopy equivalent.  Since $G$ and $H$ are minimal, it follows from
Definition~\ref{min} that they are isomorphic, and so $H$ has the
prescribed form.
\end{proof}

\section*{Acknowledgments}I would like to thank my advisor,
Sean Sather-Wagstaff, for many stimulating conversations related to
this work and for his detailed suggestions.  I would also like to
thank Lars Winther Christensen, Henrik Holm, Peter J{\o}rgensen, and
Greg Piepmeyer and for their comments and suggestions.  Finally, I
would like to thank Jan Strooker for helpful comments on the
exposition.

\bibliographystyle{amsplain}
%\bibliography{gcproj}

\providecommand{\bysame}{\leavevmode\hbox
to3em{\hrulefill}\thinspace}
\providecommand{\MR}{\relax\ifhmode\unskip\space\fi MR }
\providecommand{\MRhref}[2]{%
  \href{http://www.ams.org/mathscinet-getitem?mr=#1}{#2}
} \providecommand{\href}[2]{#2}

\end{document}